\documentclass[preprint,12pt]{elsarticle} 
\usepackage{b}
\usepackage[a4paper]{geometry}
\usepackage[small,it]{caption}

\usepackage{esint}
\usepackage{empheq}

\journal{}

\begin{document}
\begin{frontmatter}

\title{Accurate Derivative Evaluation for any Grad-Shafranov Solver}
\author[CIMS]{L.F. Ricketson}
\author[CIMS]{A.J. Cerfon}
\author[CIMS]{M. Rachh}
\author[MIT]{J.P. Freidberg}
\address[CIMS]{Courant Institute of Mathematical Sciences, New York University, New York, NY 10012}
\address[MIT]{Plasma Science and Fusion Center, Massachusetts Institute of Technology, Cambridge, MA 02139}
\date{\today}

\begin{abstract}
We present a numerical scheme that can be combined with any fixed boundary finite element based Poisson or Grad-Shafranov solver to compute the first and second partial derivatives of the solution to these equations with the same order of convergence as the solution itself. At the heart of our scheme is an efficient and accurate computation of the Dirichlet to Neumann map through the evaluation of a singular volume integral and the solution to a Fredholm integral equation of the second kind. Our numerical method is particularly useful for magnetic confinement fusion simulations, since it allows the evaluation of quantities such as the magnetic field, the parallel current density and the magnetic curvature with much higher accuracy than has been previously feasible on the affordable coarse grids that are usually implemented. 
\end{abstract}


\end{frontmatter}


\section{Introduction}

In computational physics, one often computes fields by expressing them in terms of a potential or stream function and then solving the resulting elliptic partial differential equation for the potential or stream function. The most common examples are problems involving electrostatic fields. The electric field is expressed in terms of the electric potential, which is then found via Poisson's equation. Another common situation is one in which the Euler equations for an incompressible fluid are expressed in terms of vorticity and a stream function: the desired velocity field is obtained from derivatives of the stream function, which is computed by solving Poisson's equation. A third important example, and the focus of this article, occurs in magnetic confinement fusion. The confining magnetic field in toroidally axisymmetric geometries is expressed in terms of a stream function associated with the magnetic flux, and the equilibrium magnetic configuration is computed by solving an elliptic PDE for this stream function known as the Grad-Shafranov equation \cite{grad,shafranov}.

When potentials and stream functions are computed using conventional finite difference or finite element schemes, numerical derivatives must be calculated to evaluate the fields. As a result, the convergence of the field solution is at least one order lower than the convergence order of the potential or stream function \cite{deriaz, howell}. If derivatives of the fields themselves are required, at least two orders of convergence are lost. In simulations of magnetic fusion plasmas, this loss of accuracy is a particularly acute issue. Certain physical parameters that have a key role on the stability and transport properties of hot plasmas, such as the parallel current density and the magnetic curvature, depend on first derivatives of the magnetic field or equivalently second derivatives of the flux. Popular magnetic equilibrium solvers based on a finite element formulation of the Grad-Shafranov equation compute the magnetic stream function with a numerical error that decreases as $N^{-4}$, where $N$ is the number of grid points in either direction of the two dimensional problem \cite{huysmans1991isoparametric,lutjens1992axisymmetric,lutjens1996chease,goedbloedadvanced}. The numerical errors for the derivatives of the magnetic field thus only decrease as $N^{-2}$.

There are several ways to mitigate this loss of accuracy. One can for example use higher order finite elements, as has been effectively implemented in \cite{howell}. However, in doing so, one increases the number of degrees of freedom, which leads to a larger computational cost. One can also use a spectral representation for the solution \cite{ling} or a combined modal-Green's function representation \cite{pataki2013fast}. It is indeed well known that with these methods, the loss of accuracy in the computation of derivatives is not an order of convergence, but instead a constant, as demonstrated in \citep{pataki2013fast,leeECOM} for the Grad-Shafranov equation. While satisfactory from a computational point of view, a weakness of such an approach is that the solvers rely on grids that are not typically the grids necessitated by existing stability, transport, or wave heating codes that take the computed equilibrium as an input \cite{erato,mars,jolliet,vlad}. This is the motivation for the present work. We describe a method that combines the accuracy properties of a Green's function formulation with the advantages of relying on a popular finite element solver for the computation of the solution to the Grad-Shafranov equation. The end result is the evaluation of first and second derivatives of the flux that has the same numerical accuracy and convergence properties as the flux itself while using a physically relevant grid, without having to increase the number of grid points.

Our basic idea is to obtain linear partial differential equations for the first and second partial derivatives of the potential by analytically differentiating the original PDE (say the Grad-Shafranov equation or Poisson's equation).  For instance, in the case of Poisson's equation in 2-D we see that the equations for the first derivatives are given by
\begin{equation}
	\Delta u(x,y) = f(x,y) \implies \Delta u_x = f_x(x,y), \quad \Delta u_y = f_y(x,y).
\end{equation}
Instead of computing the derivatives of $u$ from the output of the finite element solver, we use that same solver to solve the linear partial differential equations for the derivatives of $u$. 

One might then think that using the same elements as the ones used to compute $u$ automatically leads to the same order of convergence for $\nabla u$ as for $u$. This is only partially true in general because one also needs the boundary condition for the normal derivative of $u$ on the boundary of the computational domain to the same accuracy as $u$ itself. The simplest way to accomplish this task is by taking normal numerical derivatives of $u$ on the boundary. This, however, leads to the loss of one order in accuracy for each numerical derivative taken so that we are no better off than we were by simply taking numerical derivatives of the original solution over the whole domain. In other words, what is needed on the boundary is a method to compute the Dirichlet to Neumann map that does not lead to any loss in the order of convergence of the solution. The development of such a procedure is the main new contribution of the present work.

Specifically, we construct a numerical scheme that achieves the desired goal as follows. First, we re-express the Grad-Shafranov equation as a semi-linear Poisson equation with source function $F(\mathbf{x},u)$. We then decompose the solution to Poisson's equation as the sum of a particular solution $u^p$ that does not satisfy the proper Dirichlet boundary condition in general, plus a homogeneous solution $u^h$ that solves Laplace's equation and is chosen so that the full solution satisfies the proper boundary condition. We write $u^p$ as the volume integral of the product of $F(\mathbf{x},u)$ with the free-space Green's function for Poisson's equation in two dimensions. By differentiating under the integral sign, we have an exact expression for the normal derivative of $u^p$ which may be evaluated using a high-order quadrature. To compute the normal derivatives $u^h$ without ever evaluating derivatives normal to the boundary numerically, we introduce the harmonic conjugate $U$ of $u^h$ and use Green's second identity to derive a Fredholm integral equation of the second kind for $U$ on the boundary of the computational domain. After solving this integral equation for $U$ to the same order of accuracy as $u$, we can compute spectrally accurate tangential derivatives of $U$ on the boundary. Since $U$ is the conjugate gradient of $u^h$, this is equivalent to computing normal derivatives of $u^h$ with spectral accuracy, which is precisely what is needed. Note that while this work is mainly motivated by magnetic fusion applications and aimed at improving the accuracy of finite element based Grad-Shafranov solvers, our numerical method relies on a reformulation of the Grad-Shafranov equation as a Poisson problem, and can therefore also be implemented in combination with Poisson solvers.

Importantly, our approach is independent of the particular finite element solver and grid used in solving the original elliptic PDE.  However, the use of the harmonic conjugate means the approach is limited to two-dimensional problems.  The method also requires that accurate representations of the derivatives of $F(\mathbf{x},u)$ be available, and that the boundary of the computational domain be smooth. These constraints are satisfied for a large number of physically relevant problems. The extension to 3-D remains a topic for future research.

The structure of the article is as follows. In Section \ref{sec:setup}, we describe in more detail the general philosophy guiding our numerical method for obtaining the same order of convergence for the partial derivatives of $u$ as for $u$ itself. In Section \ref{sec:DtoN}, we present our scheme for the accurate computation of the Dirichlet to Neumann map. In Section \ref{sec:algorithm}, the key elements of the algorithmic implementation of our method are discussed. We demonstrate the benefits of using our scheme in combination with an existing finite element-based Grad-Shafranov solver in Section \ref{sec:examples}, in which we focus on two magnetic confinement fusion relevant examples. Finally, we summarize our work and suggest ideas for further development in Section \ref{sec:conclusion}.

\section{Setup and General Idea}\label{sec:setup}
The Grad-Shafranov equation is given by \cite{grad,shafranov}
\begin{equation} \label{gradshaf}
	r \frac{\partial}{\partial r} \left( \frac{1}{r} \frac{\partial \psi}{\partial r} \right) + \frac{\partial^2 \psi}{\partial z^2} = -\mu_0 r^2  \frac{dp}{d\psi} - \frac{1}{2} \frac{dI^2}{d\psi},
\end{equation}
with the components of the magnetic field in turn expressed as
\begin{equation} \label{psi2B}
	B_r = -\frac{1}{r} \frac{\partial \psi}{\partial z}, \quad B_\varphi = \frac{I(\psi)}{r} , \quad B_z = \frac{1}{r} \frac{\partial \psi}{\partial r}.
\end{equation}
In general, the right-hand side of Eq.\eqref{gradshaf} is a nonlinear function of $\psi$, so that the equation has to be solved iteratively \cite{pataki2013fast}. Without any penalty in terms of computational cost, one may simplify the left-hand side by the change of variables
\begin{equation}\label{transfo}
	u \coloneqq \frac{\psi}{\sqrt{r}}.
\end{equation}
This converts the Grad-Shafranov equation into the following semi-linear Poisson equation:
\begin{equation} \label{semilinPoisson}
	\Delta u = -\mu_0 r\frac{dp(\sqrt{r}u)}{du} - \frac{1}{2r} \frac{dI^2(\sqrt{r}u)}{du} + \frac{3}{4}\frac{u}{r^2}\coloneqq F(\mathbf{x},u)
\end{equation}
where $\Delta$ denotes the cartesian-coordinate Laplacian - i.e.\ 
\begin{equation}
	\Delta = \frac{\partial^2}{\partial r^2} + \frac{\partial^2}{\partial z^2}
\end{equation}
and $\mathbf{x} \coloneqq (r,z) \in \R^2$. It is worth noting that the transformation \eqref{transfo} is only advantageous in regions bounded away from the $r=0$ axis. Fortunately, this requirement holds uniformly for domains corresponding to a large class of magnetic confinement fusion applications, tokamaks in particular. 

We wish to solve equation \eqref{semilinPoisson} in a smooth, bounded domain $D$, subject to homogeneous Dirichlet boundary conditions - i.e.\ $u=0$ on $\partial D$ (which is equivalent to $\psi = 0$ on $\partial D$).  To emphasize the generality of our approach, we de-emphasize the fusion specific variables and accordingly relabel the coordinates from $(r,z) \rightarrow (x,y)$. A typical approach to solving \eqref{semilinPoisson} is to use fixed-point iteration, computing $u^n$ by
\begin{equation} \label{PoissIt}
	\Delta u^n = F \left(\mathbf{x},u^{n-1} \right),
\end{equation}
with some specified initial guess $u^0$.  One stops the iteration when
\begin{equation} \label{stop}
	\norm{ \Delta u^n - F \left( \mathbf{x}, u^n \right) } \leq \varepsilon,
\end{equation}
for some reasonable norm and specified error tolerance $\varepsilon$.  When $F(\mathbf{x},u)$ is proportional to $u$, a modified version of
the well-known inverse iteration method
\cite{trefethen1} may be used to avoid the trivial solution $u\coloneqq 0$ \cite{pataki2013fast,jardin}. Regardless of the iterative procedure, at each step, the Poisson equation (\ref{PoissIt}) has to be solved for $u^n$. While there are many techniques to solve (\ref{PoissIt}) the derivatives procedure introduced in this paper is independent of the specific technique chosen. Stated differently our aim is to present a method for the accurate computation of derivatives that can be used with \textit{any} existing finite element-based Poisson or Grad-Shafranov code. Even so, in our numerical tests, we must choose a particular technique in order to obtain results.  We use the isoparametric, bicubic Hermite finite element formulation used for the Grad-Shafranov equation in, for example, \cite{huysmans1991isoparametric,lutjens1992axisymmetric,lutjens1996chease,belien2002finesse}.  The method gives fourth order accuracy in $u$ (and thus $\psi$) and, as expected, third and second order accuracy in the first and second derivatives, respectively. We shall describe the chosen FEM scheme in sufficient detail in section 4 to understand the structure of the algorithm and the numerical results.  We refer the reader interested in more details to the sources referenced above.  

The fundamental idea behind our technique is very simple.  Say we wish to compute $u_x$.  We differentiate (\ref{semilinPoisson}) to get
\begin{equation} \label{derivPoiss}
	\Delta u_x - F_u(\mathbf{x},u) u_x = F_x(\mathbf{x},u).
\end{equation}
If the original numerical solution $u$ is known, then this is a linear elliptic PDE for $u_x$, to which the same finite element method may be applied. In this way we in principle achieve the same order of accuracy as we did for $u$. This process may be repeated by further differentiation of the starting equation, at which point another linear elliptic PDE for the second derivative may be solved to the same order accuracy.  Moreover, the original equation (\ref{semilinPoisson}) was nonlinear, and thus required iteration over many FEM solves.  In contrast, the equations for the derivatives are linear, requiring only one FEM solve.  The computation of the derivatives is thus a small additional burden relative to the original computation.

The one critical and computationally challenging point missing from the above discussion is boundary data. In order to solve (\ref{derivPoiss}), we need advance knowledge of $u_x$ on $\partial D$.  Of course, knowing $u$ to fourth order, we can easily use finite differences to compute $u_x$ on the boundary to third order, but this will lead to third order errors everywhere in the domain. We are back where we started - we know $u$ to order $k$, but can only find $u_x$ to order $k-1$. However, there is a key difference now in that the function we need to differentiate is restricted to the boundary. Note that on the boundary, $u$, $u_x$, or any higher derivatives, are smooth, periodic functions of a single surface variable.  It is thus an ideal target for \textit{spectral} differentiation.  If we can formulate the problem so that the only numerical derivatives required are spectral derivatives, we can find $u_x$ on the boundary to the same order of accuracy as $u$.  Admittedly, one may lose a constant scale factor in the accuracy of $u_x$ relative to the accuracy of $u$, since spectral differentiation inherently increases the size of high frequency modes, but this is far preferable compared to a decrease in order of accuracy.  

In the following section, we demonstrate that through the use of integral equations and a clever change of variables, it is indeed possible to compute any derivatives of $u$ on $\partial D$ to the same order of accuracy as $u$ itself while taking only spectral derivatives of periodic functions.  

\section{Boundary Data for Derivatives}\label{sec:DtoN}
To begin, it is convenient to decompose $\nabla u$ on the boundary into components $u_n$ and $u_t$, which are normal and tangent to $\partial D$, respectively.  Clearly, knowledge of the shape of $\partial D$ allows $u_x$ and $u_y$ to be computed given $u_n$ and $u_t$.  Throughout the article, we will use as our convention the \textit{inward} normal direction and the \textit{counter-clockwise} tangential direction unless otherwise noted.  

The tangential derivative is straightforward.  Direct spectral differentiation of $u$ restricted to the boundary yields $u_t$ to the same order of accuracy as $u$.  In our case, of course, $u_t \equiv 0$ is trivial to compute, but this will not be the case for second derivatives. The challenge lies in computing $u_n$. Intuitively, the difficulty comes from the fact that $u_n$ is sensitive to the behavior of $u$ \textit{off} the boundary, while $u_t$ is not.  This intuition is manifested mathematically by the fact that Green's second identity gives
\begin{equation} \label{1stkind}
	\int_{\partial D} G(\mathbf{x},\mathbf{x}') u_n(\mathbf{x}') \, dl' = \frac{1}{2} u(\mathbf{x}) - \int_D G(\mathbf{x},\mathbf{x}') F(\mathbf{x}, u(\mathbf{x}')) \, d\mathbf{x}'
\end{equation}
for any point $\mathbf{x} \in \partial D$, where $G(\mathbf{x},\mathbf{x}') = \log(\norm{ \mathbf{x} - \mathbf{x}' } )/2\pi$ is the free space Green's function of the Laplace operator.  While this does give an integral equation for $u_n$ in terms of known quantities, as desired, it is an integral equation of the first kind - the unknown only appears inside an integral - which is ill-conditioned.  Attempting to solve equation (\ref{1stkind}) for $u_n$ numerically is thus ill-advised at best, and in practice leads to very poor accuracy.  

What is needed is a transformation that converts normal derivatives to tangential derivatives, since the tangential derivatives are straightforward to compute spectrally.  For any harmonic function $\phi$ (i.e.\ $\Delta \phi = 0$), there is just such a transformation: the harmonic conjugate of $\phi$, which we denote by $\Phi$.  It is defined by
\begin{equation}
	\nabla^\perp \Phi = \nabla \phi,
\end{equation}
where $\nabla^\perp = (-\partial_y, \partial_x)^T$.  Clearly, $\phi_n = \Phi_t$, $\phi_{t}=-\Phi_{n}$ and $\Delta\phi = \Delta\Phi=0$.  The fact that $\phi$ is harmonic is essential to this transformation: if $\Delta \phi \neq 0$, then the mixed partials of $\Phi$ cannot be equal, indicating the nonexistence of such a $\Phi$.  

Unfortunately, $u$ is not harmonic.  However, it is easily decomposed into harmonic and anharmonic parts.  We write $u = u^p + u^h$, where $u^p$ is given by
\begin{equation} \label{updef}
	u^p(\mathbf{x}) = \int_D G(\mathbf{x},\mathbf{x}') F(\mathbf{x}',u(\mathbf{x}')) \, d\mathbf{x}',
\end{equation}
Here again $G(\mathbf{x},\mathbf{x}') = \log(\norm{ \mathbf{x} - \mathbf{x}' } )/2\pi$ is the free-space Green's function of the Laplace operator.  We often abbreviate this $G$ when the arguments can be understood.  The function $u^p$ is the anharmonic contribution to the total $u$. Meanwhile, the harmonic contribution $u^h$ corresponds to the homogeneous solution and
satisfies
\begin{equation}\label{BC}
	\Delta u^h = 0, \qquad \left. u^h \right|_{\partial D} = \left.-u^p\right|_{\partial D}.
\end{equation}
We can now separate the problem into the computation $u^p_n$ and $u^h_n$.  Importantly, $u^h$ is harmonic, so it has a harmonic conjugate $U^h$ satisfying $U^h_t = u^h_n$.  The strategy is thus to compute $U^h$ on $\partial D$ and to take its spectral tangential derivative to find $u^h_n$.  To find $u^p_n$, we differentiate (\ref{updef}) analytically and evaluate the resulting integral numerically. Note finally that Eq. \eqref{BC} is relatively simple because $u$ has homogeneous boundary data.  In general, the boundary condition on $u^h$ would be $\left. u^h \right|_{\partial D} = \left. u \right|_{\partial D} - \left.u^p\right|_{\partial D}$.  This general form will be necessary when computing higher order derivatives. 

The following two subsections lay out the details of computing $u^p_n$ and $u^h_n$.  Since evaluating $u^p_t$ - which can be done concurrently with evaluation of $u^p_n$ - turns out to be a key ingredient in finding $u^h_n$, we first describe the procedure for $u^p_n$ first.  In all that follows, the goal is to find values of $u^p_n$ and $u^h_n$ at points $\mathbf{x}_j \in \partial D$ that are evenly spaced in the arc-length variable $s$.  The details of computing the arc-length grid are discussed in appendix A.  

\subsection{Computing $u^p_n$}
As just mentioned, we will eventually need $u^p_t$ as well as $u^p_n$. We evaluate both from the exact formula
\begin{equation}\label{grad_part}
	\nabla u^p = \int_D \nabla G \, F(\mathbf{x}',u(\mathbf{x}')) \, d\mathbf{x}',
\end{equation}
where 
\begin{equation}
	\nabla G(\mathbf{x},\mathbf{x}') = \frac{1}{2\pi} \frac{\mathbf{x} - \mathbf{x}'}{\norm{ \mathbf{x} - \mathbf{x}' }^2}.
\end{equation}
A major advantage of the Green's function formulation is that we have an exact integral representation for the partial derivatives of $u^p$ 
given by Eq. \ref{grad_part}. This formulation, however, comes with two well known computational challenges: 
first, a robust high order quadrature rule must be used to evaluate the singular integrals; second, 
the quadrature can be potentially costly in terms of computational time since the integrals are over a two-dimensional domain. Our numerical method addresses both challenges by relying on a variant of the Quadrature By Expansion (QBX) quadrature scheme for 
the singular integrals \cite{kloeckner} accelerated by the Fast Multipole Method (FMM) \cite{greengard}. 

\begin{figure}[h!]
\centering
\includegraphics[scale=0.3]{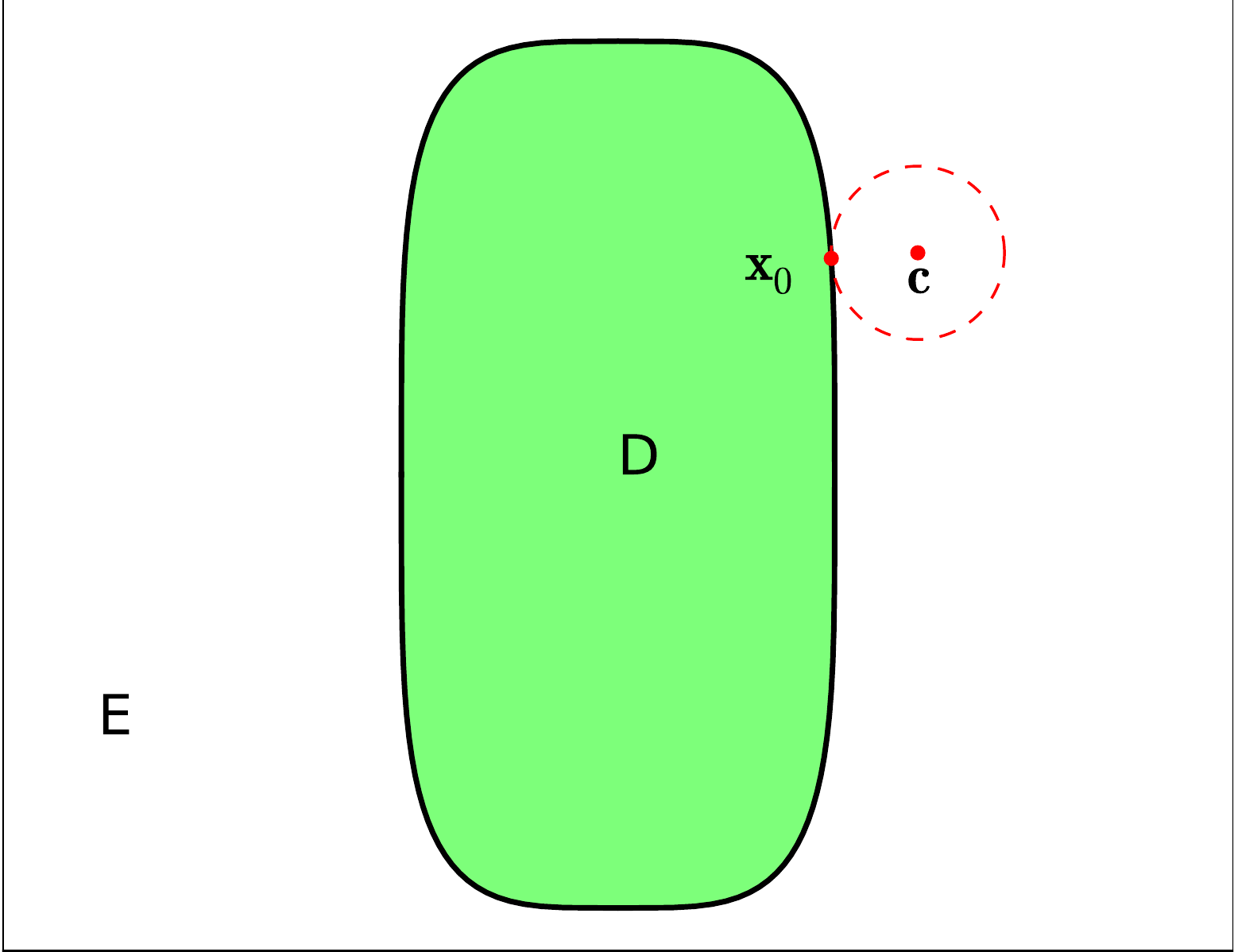}
\caption{Domain $D$ and its exterior $E$. $\mathbf{c}$ is an {\em expansion center}
corresponding to $\mathbf{x}_{0}$ in the exterior.
\label{figexpc}}
\end{figure}

Specifically, let us fix our attention to evaluating $\nabla u_{p}$ at $\mathbf{x}_{0}\in \partial D$. Let
$E=\mathbb{R}^{2} \setminus D$ be the exterior of $D$. We observe that both the components of $\nabla u_{p}$ are harmonic in $E$ and hence define smooth functions there. Using a smooth quadrature
rule, we can evaluate $\nabla u_{p}$ and its derivatives accurately at $\mathbf{c} = \mathbf{x}_{0} + r \mathbf{n}_{0}$, with $r$ a constant such that $r=O(h)$ and $h$ is
the local spacing of discretization points on the boundary (see Figure \ref{figexpc}). To make this more precise, let $\mathbf{x} = \left(x,y\right)$, $\mathbf{x}' = \left(x',y'\right)$, $\mathbf{c} = \left(
c_{x}, c_{y}\right)$, $z=x+iy$, $z'=x'+iy'$ and $c=c_{x} + i c_{y}$. We observe that $w = \partial_{x} u^{p} - i\partial_{y} u^{p}$ is complex analytic in $E$, since 
the real and imaginary parts are harmonic and complex conjugates of each other. Using equation \eqref{grad_part}, we get,
\begin{equation}
 w = \frac{1}{2\pi}\int_{D}  \frac{F(z',u(z'))}{z-z'} \, d\mathbf{x}' \, .
\end{equation}
Since $w$ is complex analytic in $E$, we can form a Taylor expansion about $c$, to obtain the following representation for $w$,
\begin{equation}
 w = \sum_{j=0}^{\infty} \left(z - c\right)^{j} \int_{D} \frac{F(z',u(z'))}{\left(z'-c\right)^{j}} \, d\mathbf{x}' \, .
\end{equation}
It can be shown that a $p$th order truncated expansion of the above equation, given by
\begin{equation}
 \tilde{w} = \sum_{j=0}^{p} \left(z - c\right)^{j} \int_{D} \frac{F(z',u(z'))}{\left(z'-c\right)^{j}} \, d\mathbf{x}' \, ,
\end{equation}
is a high order approximation to $w$, even at $\mathbf{x}_{0}$. More precisely, 
\begin{equation}
 \left|\tilde{w}\left(\mathbf{x}_{0} \right) - w \left(\mathbf{x}_{0} \right) \right| \leq C\left(p,D\right) r^{p} \, .
\end{equation}
We refer the reader to \cite{kloeckner,epstein} for a detailed discussion. To evaluate $\tilde{w}$, we still need to compute the integrals,
\begin{equation}
 \int_{D} \frac{F(z',u(z'))}{\left(z'-c\right)^{j}} \, d\mathbf{x}' \, , \label{eq:moments}
\end{equation}
which can be done to high order using a smooth quadrature rule, since the integrand is now smooth.  We choose 4th order tensor-product Gauss-Legendre quadrature on each interior grid cell - and 8th order on each boundary cell - in mapped polar-like coordinates (see (\ref{coords})).  

The computational cost of a naive implementation of the above scheme
would be $O\left(N_{b} \cdot N_{v} \right)$, where $N_{b}$ is the number of points on the boundary where we wish to evaluate $\nabla u^{p}$
and $N_{v}$ is the number of volume points used to discretize the integrands in equation \eqref{eq:moments}. However, the above
computation can be accelerated by the FMM as discussed in \cite{rachh1,rachh2} to reduce the computational cost to $O\left(N_{b} + N_{v} \right)$, which is of comparable complexity to the other steps in our method.

\subsection{Computing $u^h_n$}
Recall that to evaluate $u_{n}^h$, our approach is to compute $U^h$, the harmonic conjugate of $u^h$.  Since $U^h$ is harmonic by construction, applying Green's second identity to it gives
\begin{equation}
	\frac{1}{2} U^h(\bld{x}) = \int_{\partial D} \left( G U^h_n(\bld{x}') - G_n U^h(\bld{x}') \right) \, dl'.
\end{equation}
By rearranging and noting that $U^h_n = u^h_t = -u^p_t$ on the boundary, we have
\begin{equation} \label{seckindint}
	\frac{1}{2}U^h(\bld{x}) + \int_{\partial D} G_n U^h(\bld{x}') \, dl' = - \int_{\partial D} G u^p_t(\bld{x}') \, dl'
\end{equation}
for any $\mathbf{x} \in \partial D$.  

Given the computation of $\nabla u^p$ according to \eqref{grad_part}, it is easy to evaluate $u^p_t$, so that we can regard (\ref{seckindint}) as a second-kind integral equation for the unknown $U^h$.  We reiterate that it is crucial that this integral equation is of the second kind, in contrast to (\ref{1stkind}).  This integral equation can thus be accurately solved by discretizing each of the integrals on our arc-length grid.  It is worth mentioning that the integrand on the left-hand side of \eqref{seckindint} is not singular in spite of its appearance. This is because 
\begin{equation}
	\lim_{\bld{x'} \in \partial D \rightarrow \bld{x}} G_n(\bld{x},\bld{x}') = \frac{1}{4\pi} \kappa (\bld{x}),
\end{equation}
where $\kappa$ is the signed curvature of $\partial D$.  Given a parameterization of the curve $(x(t),y(t))$, it is defined by
\begin {equation}
	\kappa \coloneqq \frac{\ddot{x} \dot{y} - \ddot{y} \dot{x}}{(\dot{x}^2 + \dot{y}^2)^{3/2}}.
\end{equation}

Thus, $G_n$ restricted to the boundary is a smooth function once we define $G_n(\mathbf{x},\mathbf{x}) \coloneqq \kappa(\mathbf{x})/4\pi$.  As such, we are free to define
\begin{equation}
	\beta_{ij} \coloneqq \Delta s \, G_n(\mathbf{x}_i,\mathbf{x}_j),
\end{equation}
where $\Delta s$ is the arc-length distance between grid points.  Thus, the approximation of the integral on the left-hand side of \eqref{seckindint} by the trapezoidal rule is written as
\begin{equation}
	\int_{\partial D} G_n U^h(\mathbf{x}') \, dl' \approx \sum_j \beta_{ij} U^h_j.
\end{equation}
The trapezoidal rule suffices since it is well known to converge spectrally for smooth, \textit{periodic} functions \cite{trefethen2}.  

The integrand on the right-hand side of (\ref{seckindint}), on the other hand, is logarithmically singular.  We denote by $\gamma_j$ the value of the integral 
at $\mathbf{x} = \mathbf{x}_j$, and again compute it using the standard version of QBX described above.  

In the end, the discretized version of (\ref{seckindint}) is
\begin{equation} \label{seckinddisc}
	\sum_j \left( \frac{1}{2} \delta_{ij} + \beta_{ij} \right) U^h_j = \gamma_i.
\end{equation}
This linear system is ill-conditioned because $U^h$ is defined in terms of its derivatives, and thus only well defined up to a constant. We thus expect infinitely many solutions to (\ref{seckinddisc}). The problem is solved by imposing the additional constraint that $U^h$ average to zero. That is,
\begin{equation} \label{zeroavg}
	\Delta s \sum_j U^h_j = 0.
\end{equation}
We may impose this constraint without increasing the size of the system by solving
\begin{equation} \label{seckindintmod}
	\sum_j \left( \frac{1}{2} \delta_{ij} + \beta_{ij} + \Delta s^2 \right) U^h_j = \gamma_i.
\end{equation}
As shown in \cite{sifuentes2014randomized}, this new linear system has a unique solutions which also solves Eq.\ (\ref{seckindint}) with probability 1.   

The linear system (\ref{seckindintmod}) may now be solved using any standard linear algebra package.  It bears mentioning that this system is dense, while the stiffness matrix in the FEM formulation is sparse.  It may thus seem that solving (\ref{seckindintmod}) becomes the rate limiting step in the computation.  However, this is not the case.  The dimension of the dense system scales with the number of grid points on the boundary, while the stiffness matrix size scales with the number of points in the entire two dimensional domain - roughly the square of the number of boundary points.  Thus, the complexity of solving each of the two linear systems is comparable.  Moreover, the system (\ref{seckindintmod}) need only be solved once, while the finite element system must in general be solved several times since the semi-linear Poisson equation \eqref{semilinPoisson} has to be solved iteratively.  

\section{Implementation Details and Algorithm Summary}\label{sec:algorithm}
As mentioned in Section \ref{sec:setup}, we choose an isoparametric, bicubic Hermite finite element formulation, which we briefly outline here.  The computational domain is represented in $(\rho,\theta)$ coordinates, which are defined in terms of the Cartesian coordinates $x$ and $y$ by the following transformation
\begin{equation} \label{coords}
\begin{split}
	x &= x_c + \rho f(\theta) \cos \theta \\
	y &= y_c + \rho f(\theta) \sin \theta
\end{split}
\end{equation}
for some pre-specified central axis $(x_c,y_c)$. Accordingly, $\theta$ is the polar angle, and $f(\theta)$ is represented by a Fourier series, which is computed from given boundary information.  When necessary, tangent and normal directions to the boundary may be found by differentiating this Fourier series.  

In $(\rho,\theta)$ coordinates, the computational domain is $[0,1] \times [0,2\pi)$.  This is subdivided into an $N\times N$ grid, to which the FEM is applied using the bicubic Hermite basis in the $(\rho,\theta)$ variables.  The integrals of the basis functions that must be evaluated are computed using tensor-product Gauss-Legendre quadrature with four nodes in each direction.  The solution $u$ of (\ref{semilinPoisson}) is approximated by iteratively applying this FEM scheme to (\ref{PoissIt}), stopping according to the criterion (\ref{stop}), using infinity norm over the unknowns and $\varepsilon = 1.7 \times 10^{-14}$.  

For the purposes of spectral differentiation, it is desirable to have a grid along the boundary that is evenly spaced in arc length.  The angles $\theta_j$ corresponding to such a grid may be computed using the method described in appendix A.  We construct this grid using $8N$ points, and denote the arc-length grid-points on the boundary by $\mathbf{x}_j$.  We proceed through the following steps to find $u_x$ and $u_y$ on the boundary.  Throughout, $\mathbf{n}_j = (n_{x,j}, n_{y,j})$ denotes the unit inward normal to the boundary at $\mathbf{x}_j$, and $\mathbf{t}_j$ the unit counterclockwise tangent vector.  We present the method in the more complex case in which it is not assumed that $u = 0$ on the boundary, for this general case is required when computing second derivatives.  
\begin{enumerate}
	\item Use FMM accelerated QBX to evaluate
	\begin{equation} \label{evalthis}
		\nabla u^p_j= \int_{D} \nabla G(\mathbf{x}_j,\mathbf{x}') F(\mathbf{x}',u(\mathbf{x}')) \, d\mathbf{x}',
	\end{equation} as described in 3.1.  
	\item Set $u^p_{n,j} = \mathbf{n}_j \cdot \nabla u^p$ and $u^p_{t,j} = \mathbf{t}_j \cdot \nabla u^p(\mathbf{x}_j)$.  
	\item Compute
	\begin{equation}
		\gamma_j = \int_{\partial D} G(\mathbf{x}_j,\mathbf{x}') \left[u_t(\mathbf{x}') - u^p_t(\mathbf{x}')\right] \, dl'
	\end{equation}
	by again using FMM accelerated QBX applied to the computed values of $u^p_{t,j}$.  Here, $u_t$ is computed by direct spectral differentiation of the FEM solution.  
	\item Solve the linear system
	\begin{equation}
		\sum_j \left( \frac{1}{2} \delta_{ij} + \beta_{ij} + \Delta s^2 \right) U^h_j = \gamma_i,
	\end{equation}
	where $\Delta s$ is the spacing of the arc-length grid, as described in 3.2. In our code, we simply use MATLAB's ``backlash" operator, as we find this is not the rate limiting step in the procedure. More generally, an iterative solver based on the generalized minimal residual method (GMRES) is satisfying since the linear system only needs to be inverted once. The system is well conditioned so will converge quickly, i.e in $O(1)$ iterations.
	\item Compute $U^h_{t,j}$ by spectral differentiation of $U^h_j$.  
	\item Set $u_{n,j} = U^h_{t,j} + u^p_{n,j}$, from which it follows that $u_{x,j} = n_{x,j} u_{n,j} + n_{y,j} u_{t,j}$, and $u_{y,j} = n_{y,j} u_{n,j} - n_{x,j} u_{t,j}$.  
\end{enumerate}
With this boundary data in hand, we use the same finite element formulation to solve
\begin{equation}
\begin{split}
	\Delta u_x - F_u(\mathbf{x},u) u_x &= F_x(\mathbf{x},u) \\
	\Delta u_y - F_u(\mathbf{x},u) u_y &= F_y(\mathbf{x},u)
\end{split}
\end{equation}
for $u_x$ and $u_y$, respectively.  

Computing second derivatives is directly analogous.  The procedure above is re-used, but with $u \rightarrow u_x$ everywhere and $F(\mathbf{x,u}) \rightarrow F_x(\mathbf{x},u) + F_u(\mathbf{x},u) \mathbf{u}_x$ in (\ref{evalthis}).  To be more specific, we evaluate
\begin{equation}
	\nabla u_x^p = \int_D \nabla G(\mathbf{x},\mathbf{x}') \left[ F_x(\mathbf{x}',u(\mathbf{x}') + F_u(\mathbf{x}',u(\mathbf{x}'))u_x(\mathbf{x}') \right] \, d\mathbf{x}'
\end{equation}
on the arc-length grid, and
solve
\begin{equation}
	\frac{1}{2}U_x^h(\bld{x}) + \int_{\partial D} G_n U_x^h(\bld{x}') \, dl' = \int_{\partial D} G(\mathbf{x},\mathbf{x}') \left[ u_{xt}(\mathbf{x}') - u^p_{xt}(\mathbf{x}') \right] \, dl'
\end{equation}
for $U_x^h$.  Then, $u_{xn} = U^h_{xt} + u^p_{xn}$ and $u_{xt}$ is found by direct spectral differentiation of $u_x$ on the same grid.  

Having $u_{xn}$ and $u_{xt}$, we may find $u_{xx}$ and $u_{xy}$ on the boundary, which are used as boundary data to solve
\begin{equation}
\begin{split}
	\Delta u_{xx} - F_u(\mathbf{x},u) u_{xx} &= F_{xx} + 2F_{xu}u_x + F_{xu} u_x^2 \\
	\Delta u_{xy} - F_u(\mathbf{x},u) u_{xy} &= F_{xy} + F_{xu} u_y + F_{yu} u_x + F_{xy} u_x u_y
\end{split}
\end{equation}
for $u_{xx}$ and $u_{xy}$ respectively.  It is then simple to find $u_{yy}$, since
\begin{equation}\label{eq:uyy}
	u_{yy} = F(\mathbf{x},u) - u_{xx}
\end{equation}
as a consequence of the original PDE (\ref{semilinPoisson}).  Once $u$ and its derivatives are known, it is straightforward to compute $\psi$ and its derivatives to the same accuracy using the formula \eqref{transfo}.  
\section{Numerical Results}\label{sec:examples}
For numerical tests, we consider two situations in which exact solutions are known. The first situation corresponds to toroidally axisymmetric plasma equilibria with pressure and current profiles chosen so that the right-hand side of the Grad-Shafranov equation does not depend on $\psi$. Exact solutions are compared with numerical solutions for two magnetic confinement devices, as discussed in detail in section \ref{sec:solovev}. In this first situation, both equations \eqref{gradshaf} and \eqref{semilinPoisson} are linear partial differential equations. In order to show that our procedure applies just as well to nonlinear partial differential equations, we consider in section \ref{sec:nonlinear} an exact solution to the nonlinear Poisson-Boltzmann equation on an ellipse-like domain.

\subsection{Grad-Shafranov equation with Solov'ev profiles}\label{sec:solovev}
In this first example, we solve the Grad-Shafranov equation with a Solov'ev pressure profile \cite{solovev, cerfon2010one} and no diamagnetic or paramagnetic contribution to the toroidal magnetic field ($dI^2/d\psi=0$). The simplified Grad-Shafranov equation is given by
\begin{equation}\label{solovev}
	x \frac{\partial}{\partial x} \left( \frac{1}{x} \frac{\partial \psi}{\partial x} \right) + \frac{\partial^2 \psi}{\partial y^2} = Cx^2
\end{equation}
A simple and exact solution to this equation is
\begin{equation}
	\psi = \frac{C}{8} x^4 + d_1 + d_2 x^2 + d_3(x^4 - 4x^2y^2)
\end{equation}
for any constants $d_1$, $d_2$, and $d_3$, where the boundary $\partial D$ corresponding to $\psi=0$ is simply chosen to be the curve along which the above polynomial vanishes. Of course, not every choice of $C$, $d_1$, $d_2$, and $d_3$ results in a reasonable plasma boundary. We choose these constants following the approach in \cite{pataki2013fast,cerfon2010one}.  Namely, we rewrite them in terms of plasma relevant quantities.  These are $\epsilon$, $\delta$, and $\kappa$, which are called the inverse aspect ratio, triangulation, and elongation, respectively.  By imposing
\begin{equation}
	\psi(1+\epsilon,0) = \psi(1-\epsilon) = \psi(1-\delta\epsilon,\kappa\epsilon) = 0,
\end{equation}
we may write an invertible linear system that relates $(\epsilon,\delta,\kappa)$ to $(d_1,d_2,d_3)$ \cite{pataki2013fast}.  In all our examples, we will use $C = 10$. Its value is not important since it can be scaled out of the problem.

We compute the errors in $u$ and its derivatives using the $L^2$ norm.  That is, the error in any quantity $Q(\mathbf{x})$ is given by
\begin{equation}
	\textrm{Error} = \left(\int_D (Q_{\textup{approx}} - Q_{\textup{exact}})^2 \, d\mathbf{x} \right)^{1/2}.
\end{equation}
The integral is evaluated using Gauss-Legendre quadrature over the FEM grid in $(s,\theta)$ coordinates.  In addition, we measure error in the location of the magnetic axis, which is the value of $\mathbf{x}$ at which $\nabla \psi = \mathbf{0}$. This is computed via gradient descent on the $y=0$ axis, since symmetry tells us the magnetic axis must lie on this line. The numerical result is compared to the exact axis location, given by 
\begin{equation}
	x^* = 2\sqrt{\frac{-d_2}{C + 8d_3}}.
\end{equation}

\begin{figure}[h!]
  \centering
	\includegraphics[width=.49\textwidth]{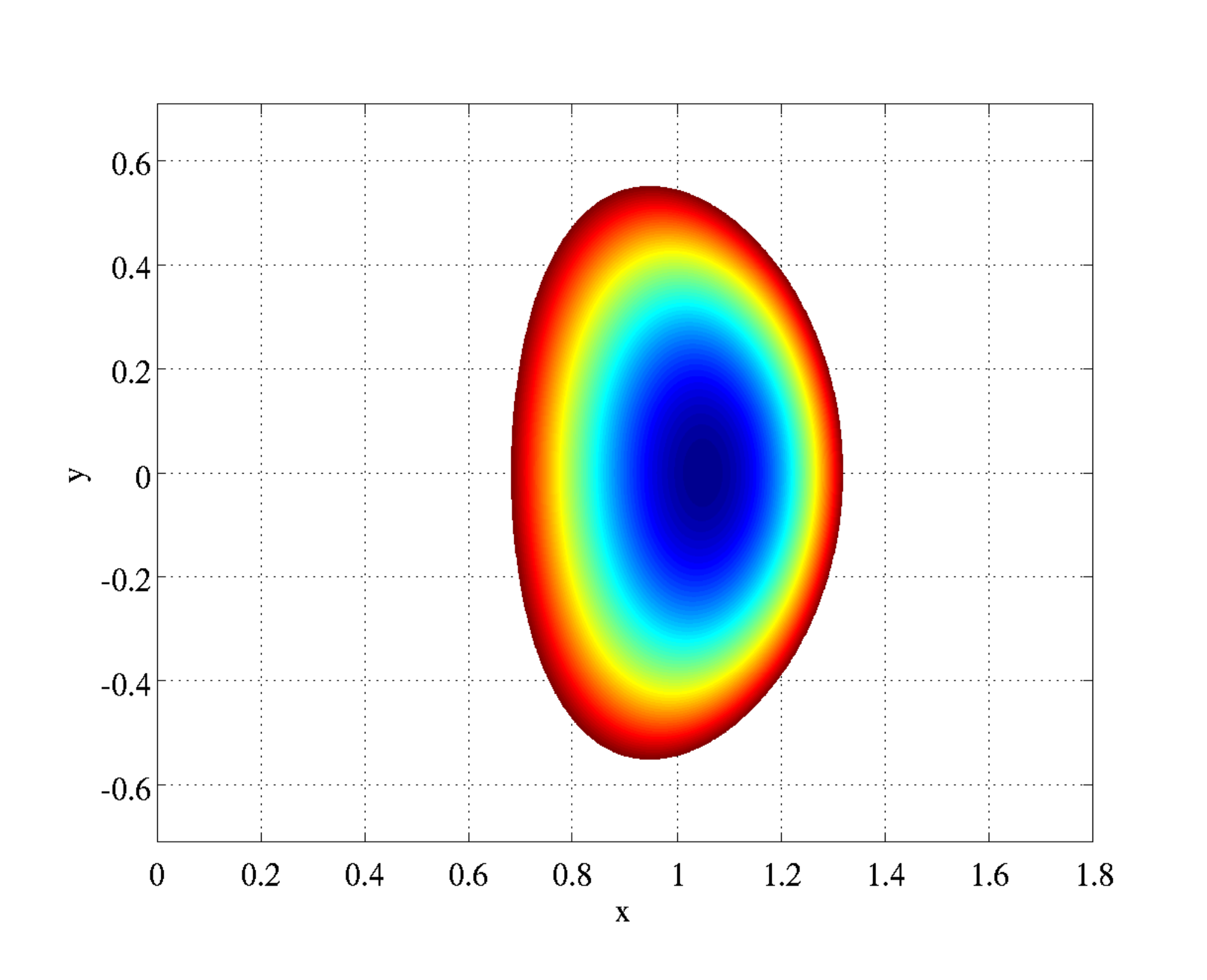}
	\includegraphics[width=.49\textwidth]{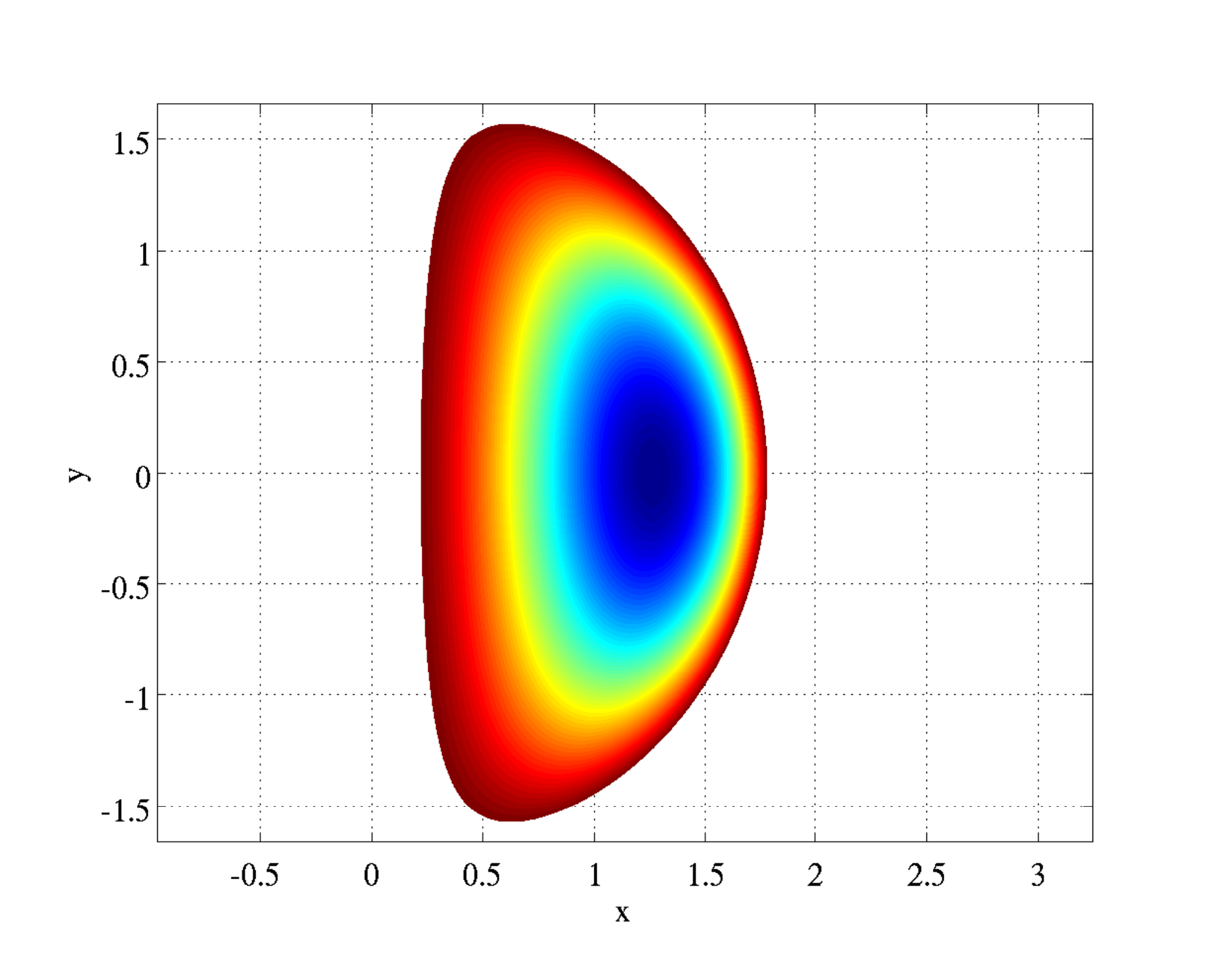}
	\caption{\textit{Left}: Contours of $\psi$ for the ITER-like Solov'ev equilibrium given by Eq.(\ref{solovev}) with $(\epsilon,\delta,\kappa) = (0.32, 0.33, 1.7)$.  \textit{Right}: Contours of $\psi$ for the NSTX-like Solov'ev equilibrium given by Eq.(\ref{solovev}) with $(\epsilon,\delta,\kappa) = (0.78, 0.35, 2)$.}
	\label{equilibria}
\end{figure}

We consider two cases: 1) $(\epsilon,\delta,\kappa) = (0.32, 0.33, 1.7)$, which corresponds to the geometry of the ITER tokamak \cite{ITER}; 2) $(\epsilon,\delta,\kappa) = (0.78, 0.35, 2)$, which corresponds to the geometry of the spherical tokamak NSTX \cite{NSTX}. For reference, the equilibria are shown in figure \ref{equilibria}.  Since the equations for the Solov'ev equilibria are in fact linear, we perform an additional test on a Poisson-Boltzmann equation to confirm the method's performance on nonlinear equations.  

\subsubsection{ITER example}\label{sec:ITER}

\begin{figure}[h!] 
  \centering
	\includegraphics[width=.7\textwidth]{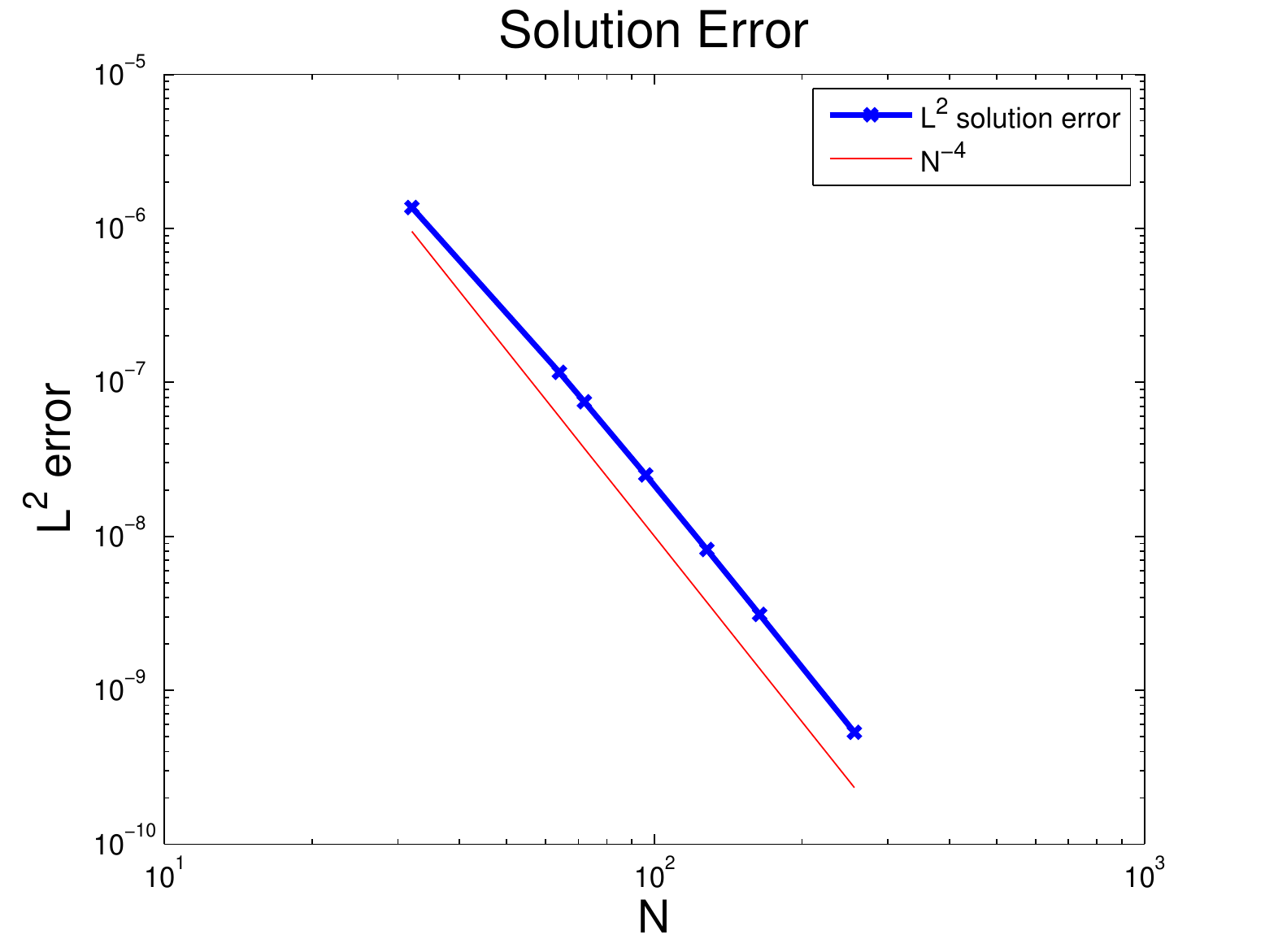} 
	\caption{Error in the solution $u$ for the ITER-like Solov'ev equilibrium given by Eq.(\ref{solovev}) with $(\epsilon,\delta,\kappa) = (0.32, 0.33, 1.7)$, displaying the expected $N^{-4}$ convergence rate.}
	\label{ITER1}
\end{figure}

\begin{figure}[h!]
  \centering
	\includegraphics[width=.49\textwidth]{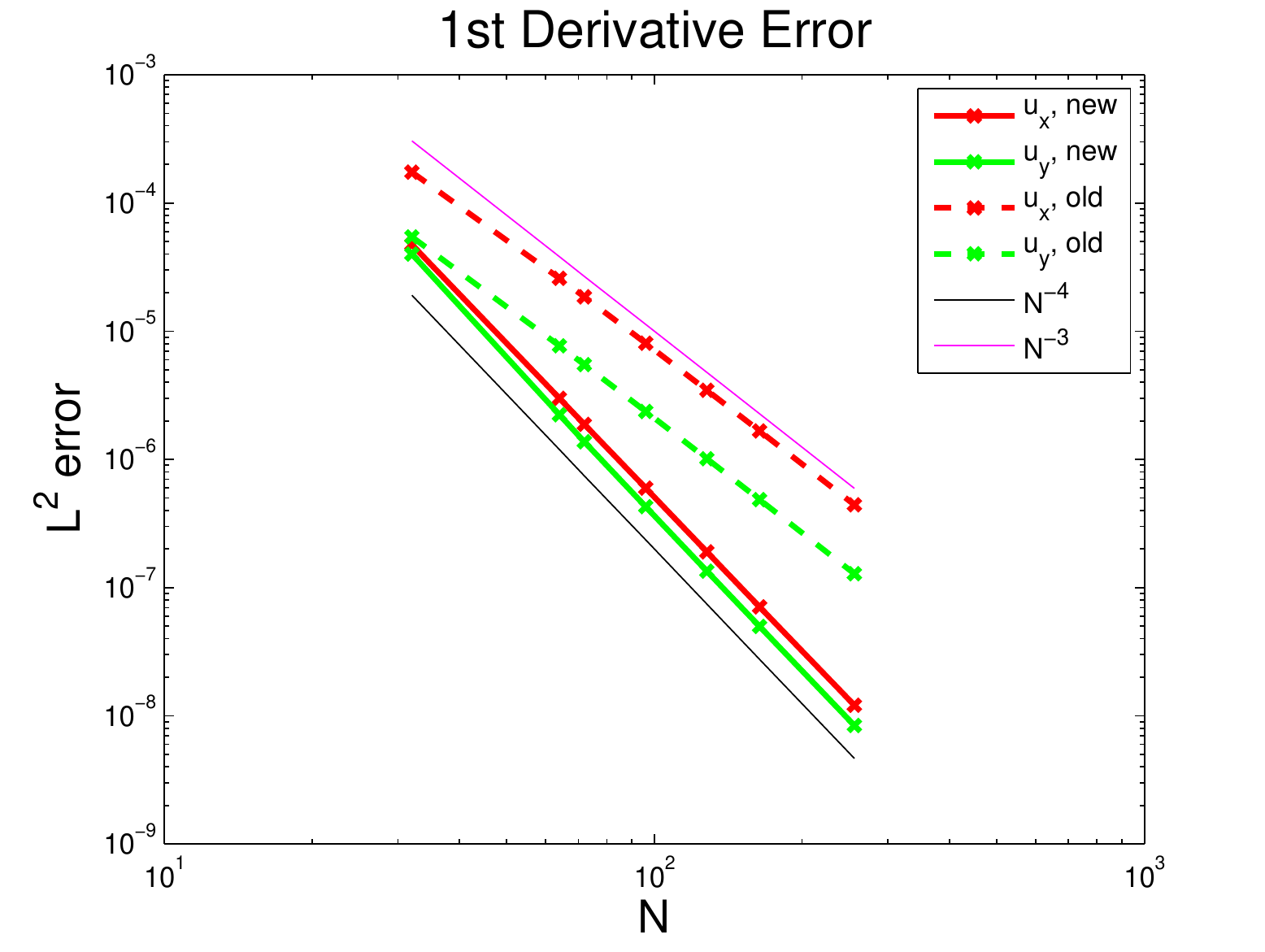}
	\includegraphics[width=.49\textwidth]{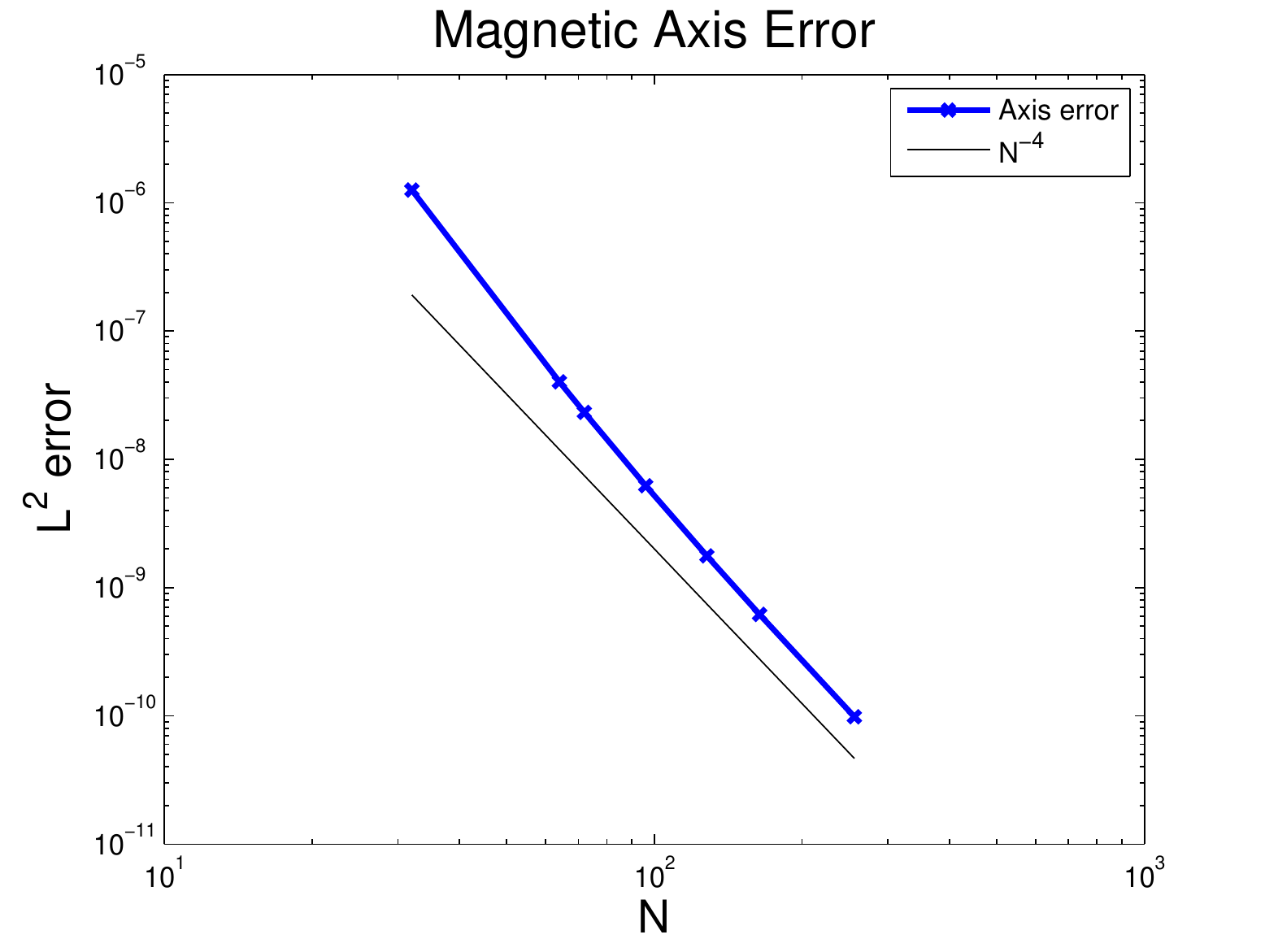}
	\caption{\textit{Left}: Error in the first partial derivatives for the ITER-like Solov'ev equilibrium given by Eq.(\ref{solovev}) with $(\epsilon,\delta,\kappa) = (0.32, 0.33, 1.7)$.  Dashed lines represent the naive computation, while results obtained with the new method presented in this article are shown in solid lines.  \textit{Right}: Error in the location of the magnetic axis, displaying improved convergence relative to results directly obtained with CHEASE and reported in \cite{lutjens1992axisymmetric,lutjens1996chease}.}
	\label{ITER2}
\end{figure}

Let us first focus on the ITER-like equilibrium. In figure \ref{ITER1}, we plot the error in the solution $u$.  As expected, we see that this error is $O(N^{-4})$ as a result of the FEM scheme used.  More interestingly, in figure \ref{ITER2} we plot the error in $u_x$ and $u_y$, as well as the location of the magnetic axis. The dashed lines are computed by reading off $u_\theta$ and $u_s$ from the finite element solution, and then converting to rectangular coordinates.  As expected, one order of accuracy is lost relative to the solution.  

In contrast, the solid lines are the errors using the new method presented in this article.  One can observe the improved convergence rate and improved absolute error for any $N > 32$.  At $N = 256$, we see an improvement of more than an order of magnitude in each derivative.  The error in the magnetic axis is also $O(N^{-4})$.  This is to be compared with the results presented in \cite{lutjens1992axisymmetric}, which did not take advantage of our method, and where $O(N^{-3})$ accuracy was observed using the same finite element basis.  

\begin{figure}[h!]
  \centering
	\includegraphics[width=.7\textwidth]{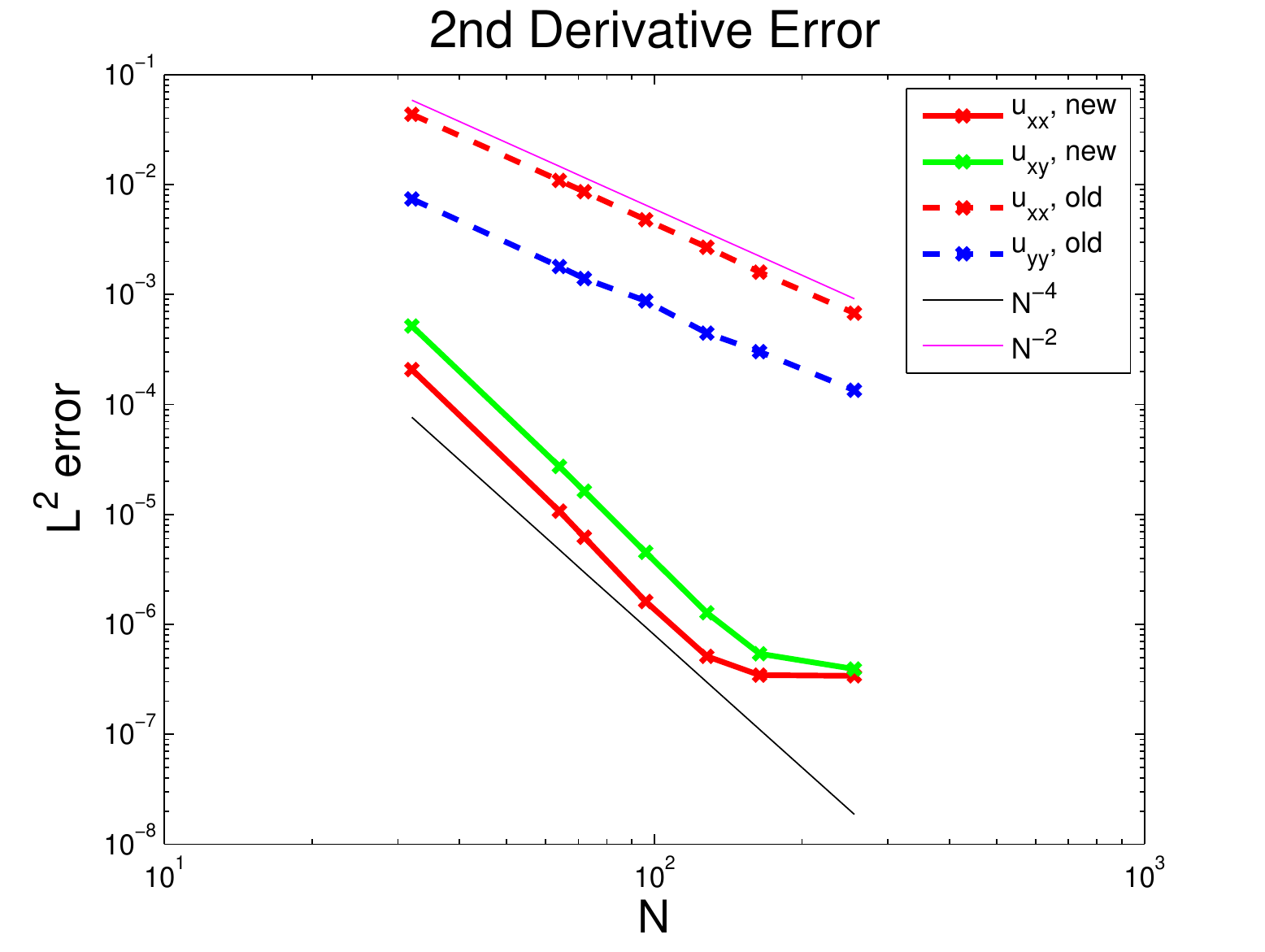}
	\caption{Error in the second derivatives for the ITER-like Solov'ev equilibrium given by Eq.(\ref{solovev}) with $(\epsilon,\delta,\kappa) = (0.32, 0.33, 1.7)$. Dashed lines represent the naive computation, while results obtained with the new method presented in this article are shown in solid lines.}
	\label{ITER3}
\end{figure}

Even more significant are the improvements in accuracy of the second derivatives, shown in figure \ref{ITER3}.  For the clarity of the figure, we did not plot the error in $u_{yy}$ obtained with the new method because it is similar to the error in $u_{xx}$ according to equation \eqref{eq:uyy}. We did not plot the error in $u_{xy}$ that one obtains with the standard finite element approximation either because we could not get it to converge properly. We note by looking at figure \ref{ITER3} that the improved convergence rate with the new method is observed for moderate $N$, but the curve flattens as $N$ grows beyond $\approx 128$. This is due to round-off error that accumulates from taking numerous spectral derivatives. Even so, at the present convergence rate, the standard finite element approximation requires $N \approx 5000$ to reach the accuracy the new method achieves at $N = 256$.  That is a reduction in complexity of roughly $(5000/256)^2 \approx 381$, which easily compensates for the additional computation resulting from execution of the methods presented here. An alternate comparison is that for $N = 128$ the new method is more than 1000 times more accurate for $u_{xx}$ than the direct FEM result.
 
Even in the event that only 4 digits of accuracy are required in the second derivatives, the new scheme accomplishes this with $N\approx 50$, while the standard method requires $N \approx 315$.  The speed improvement in this case is roughly $(315/50)^2 \approx 40$.  This again easily outstrips the additional cost devoted to the boundary and the extra finite element solves.

\subsubsection{NSTX example}\label{sec:NSTX}
\begin{figure}[h!]
  \centering
	\includegraphics[width=.49\textwidth]{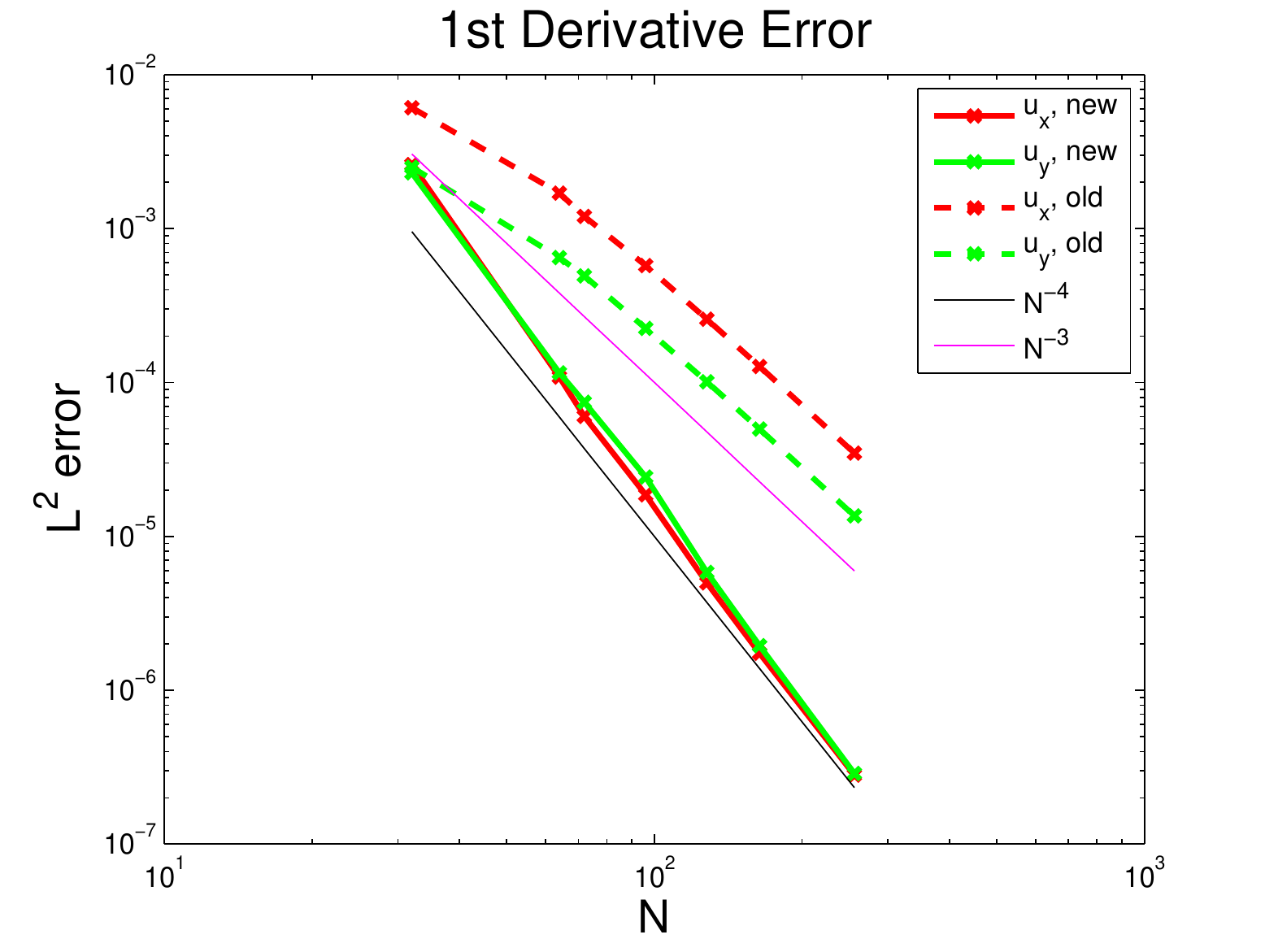}
	\includegraphics[width=.49\textwidth]{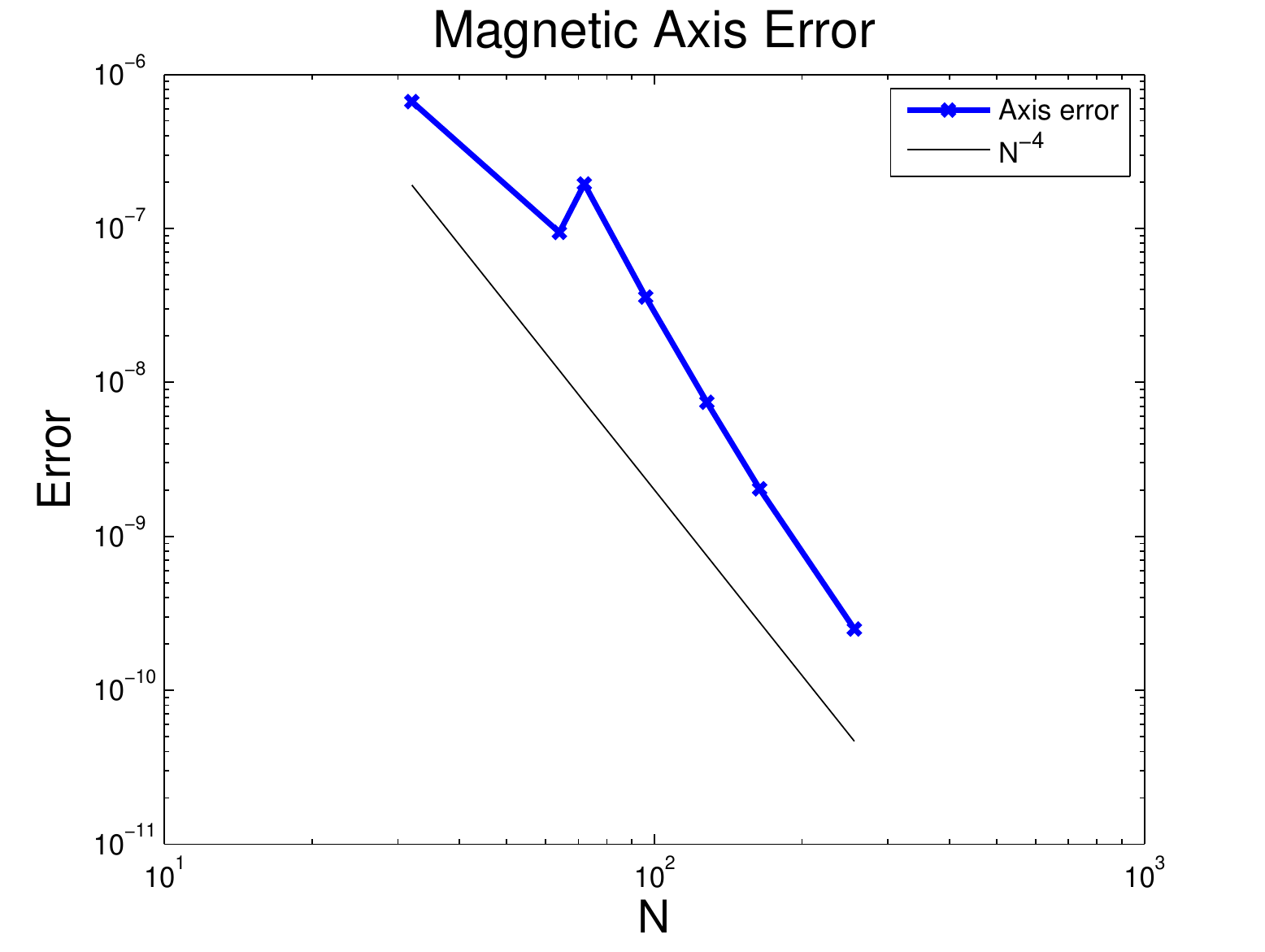}
	\caption{\textit{Left}: Error in the first partial derivatives for the NSTX-like Solov'ev equilibrium given by Eq.(\ref{solovev}) with $(\epsilon,\delta,\kappa) = (0.78, 0.35, 2)$.  Dashed lines represent the naive computation, while results obtained with the new method presented in this article are shown in solid lines.  \textit{Right}: Error in the location of the magnetic axis.}
\label{NSTX1}
\end{figure}

\begin{figure}[h!]
  \centering
	\includegraphics[width=.7\textwidth]{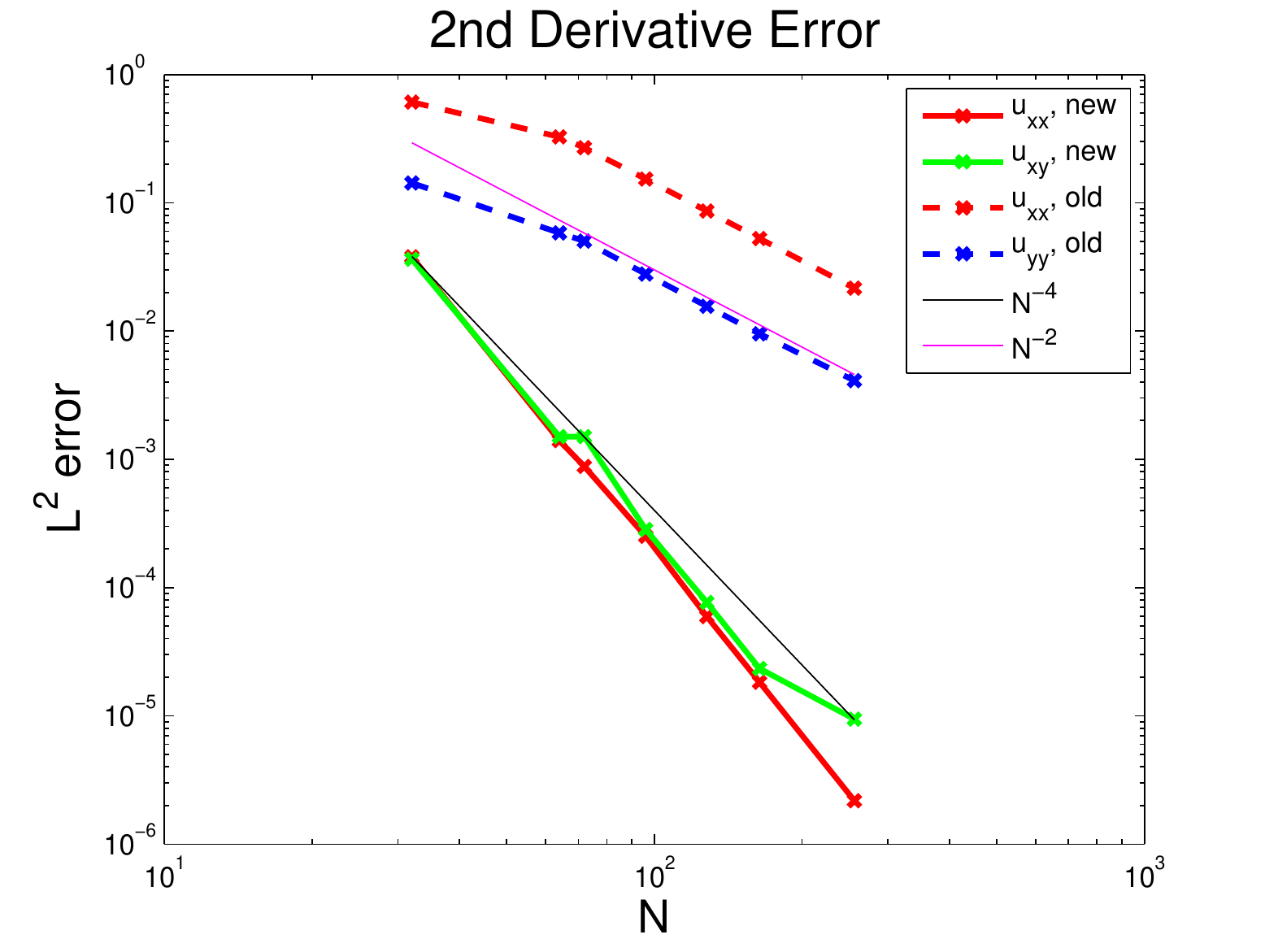}	
	\caption{Error in the second derivatives for the NSTX-like Solov'ev equilibrium given by Eq.(\ref{solovev}) with $(\epsilon,\delta,\kappa) = (0.78, 0.35, 2)$. Dashed lines represent the naive computation, while results obtained with the new method presented in this article are shown in solid lines.}
\label{NSTX2}	
\end{figure}

As a second numerical example, we consider the more challenging case of a Solov'ev equilibrium as given by Eq.\eqref{solovev} for the high inverse aspect ratio, highly elongated spherical tokamak NSTX \cite{NSTX}, with typical parameters $(\epsilon,\delta,\kappa) = (0.78, 0.35, 2)$. The results are shown in figures \ref{NSTX1} and \ref{NSTX2}. We did not plot the curves corresponding to $u_{yy}$ with the new method and $u_{xy}$ with the standard method for the same reasons as the ones we gave in section \ref{sec:ITER}. We observe the same improvement in terms of the order of convergence of the derivatives of the solution to the Grad-Shafranov equation. One does note that for the magnetic axis and second derivatives, the convergence is not as smooth at small $N$ as in the ITER case. We hypothesize that this is due to the fact that the NSTX boundary is more difficult to resolve so that larger $N$ is required to observe the asymptotic behavior.  

\subsection{Poisson-Boltzmann example}\label{sec:nonlinear}
We now verify that the method we describe in this article also performs well for nonlinear cases by considering the nonlinear Poisson problem given by
\begin{equation}
	\Delta u = \alpha e^{-u}, \qquad \left. u \right|_{\partial D} = 0,
\end{equation}
where $\partial D$ is defined by
\begin{equation}
	c_1 \cosh ky - c_2 \cos kx = 1.
\end{equation}
We pick $\alpha = 2k^2(c_1^2 - c_2^2)$, giving an exact solution,
\begin{equation}
	u = 2 \log (c_1 \cosh ky - c_2 \cos kx),
\end{equation}
against which our numerical solutions may be compared. We evaluate the numerical error in the same way as explained in section \ref{sec:solovev}. For the free parameters, we choose $k = \pi/5$, $c_1 = .0287$ and $c_2 = .3301$.  This gives an ellipse-like boundary with a 2-to-1 aspect ratio. 

\begin{figure}[h!]
  \centering
	\includegraphics[width=.49\textwidth]{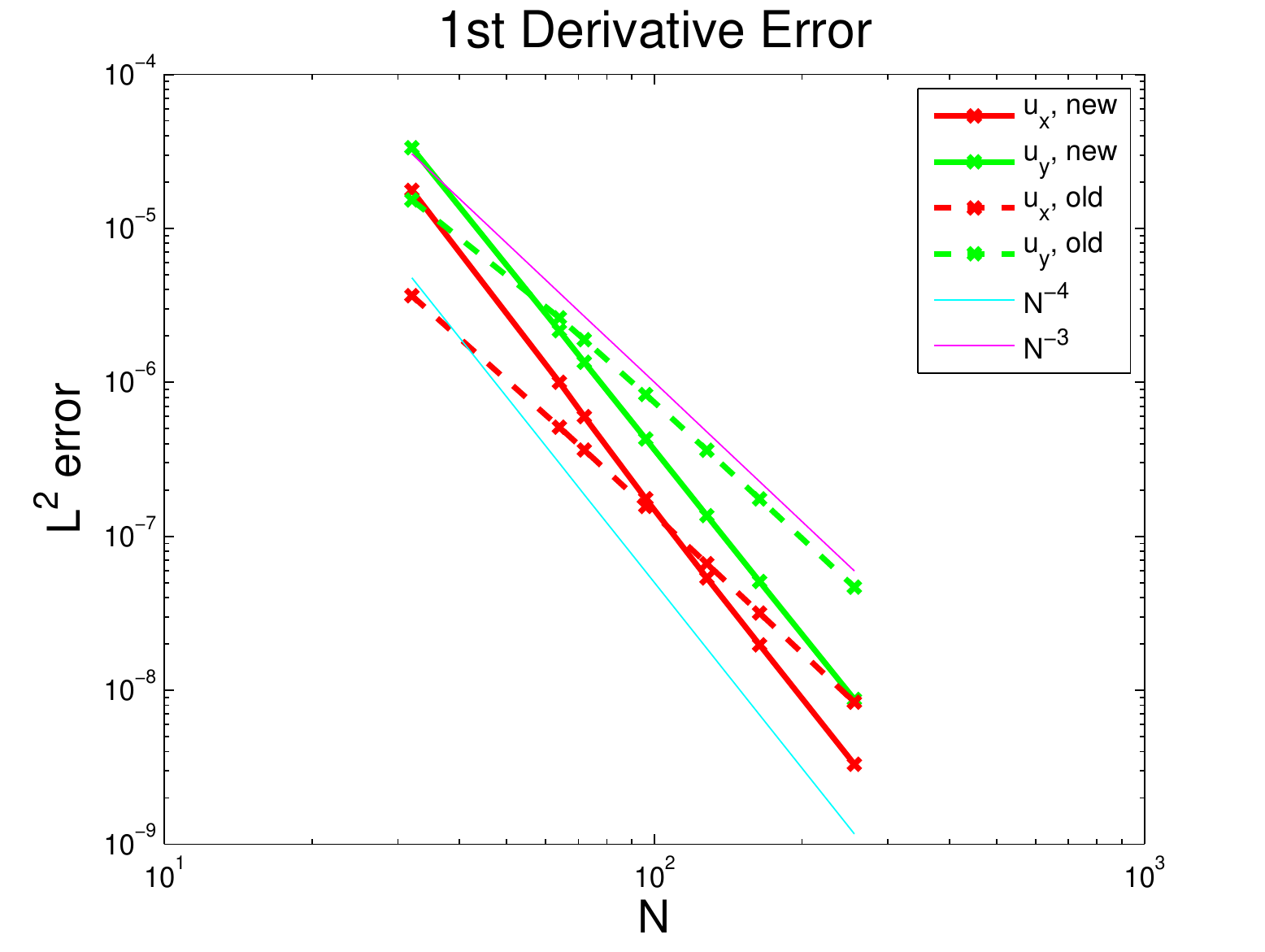}
	\includegraphics[width=.49\textwidth]{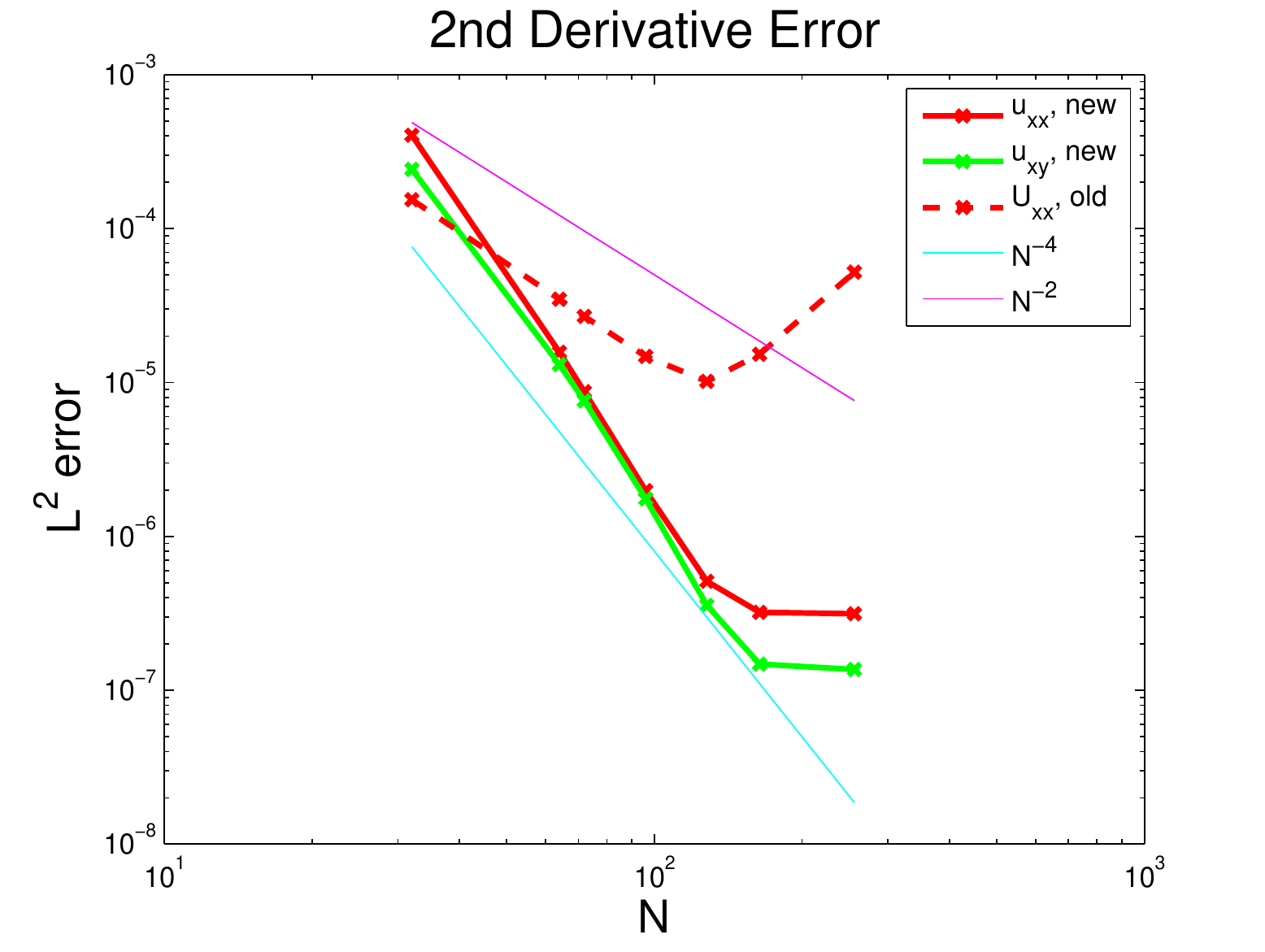}
	\caption{\textit{Left}: Error in the first partial derivatives for the Poisson-Boltzmann example.  Dashed lines represent the naive computation, while results obtained with the new method presented in this article are shown in solid lines.  \textit{Right}: Error in the second derivatives, with dashed lines again corresponding to the naive computation and solid lines to the new method.}
	\label{Boltzmann}
\end{figure}

Results for the first and second derivatives are plotted in figure \ref{Boltzmann}. The method presented here performs as expected.  We again observe the flattening of the second derivative curve due to accumulated round-off error from numerous spectral derivatives.  Curiously, the finite difference approximation of second derivatives behaves even worse than expected.  We've plotted $u_{xx}$, which initially converges at the expected rate before diverging at large $N$.  However, $u_{xy}$ and $u_{yy}$ fail to converge at all when finite differences are used.  This hints that in addition to being more accurate, the new method may prove more robust than the standard approach.  

\section{Discussion and Conclusions}\label{sec:conclusion}
We have presented a method based on an integral equation formulation for the accurate numerical evaluation of the Dirichlet to Neumann map for Poisson's equation and the Grad-Shafranov equation. The method takes as an input the solution to either of these partial differential equations obtained with any finite element-based solver. By differentiating the partial differential equation analytically, we are then able to use the same finite element solver to compute partial derivatives of the solution with the same order of convergence as the solution itself. The computational cost of implementing the numerical method we describe in this article is comparable to two Poisson or Grad-Shafranov solves. This additional computational cost is offset by the fact that on small grids, we obtain accuracies for the first and second derivatives that could only be achieved on a much larger grid with the standard method. In the context of plasma physics simulations, a code based on our numerical method can be advantageously combined with existing and commonly used Grad-Shafranov solvers to calculate quantities such as the magnetic field, the parallel current density and the magnetic curvature with much higher accuracy than is possible with the standard methods on the fairly small grids that can usually be afforded. At present, our method is only applicable to cases for which the boundary of the domain is smooth. The quadrature schemes we employ here can handle boundaries with a sharp corner, but the Fourier based method we use to compute spectrally accurate derivatives on the boundary would not yield accurate results for curves that are not smooth. Given the importance of boundaries with corners in magnetic confinement fusion, an improved scheme with this capability, most likely based on local interpolation instead of global interpolation, would be desirable. This is the subject of ongoing work, with progress to be reported at a later date.


\section*{Acknowledgments}
The authors would like to thank Prof.\ Leslie Greengard (NYU CIMS) for many insightful conversations and Dr.\ Carlos Borges (NYU CIMS) for his help in implementing the QBX codes in MATLAB. L.F.R. and A.J.C. were supported by the U.S. Department of Energy, Office of Science, Fusion Energy Sciences under Award Nos. DE-FG02-86ER53223 and DE-SC0012398. M.R. was supported in part by the Department of Energy, Office of Science, Office of Advanced Scientific Computing Research, Applied Mathematics Program under Award DEFGO288ER25053, and by the Office of the Assistant Secretary of Defense for Research and Engineering and AFOSR under NSSEFF Program Award FA9550-10-1-0180.

\begin{appendix} 


\section{Arc-length grid computation}
We desire a sequence of angles $\theta_j$ that correspond to equally spaced points in arc-length, with spacing given by $\Delta s = 8N / L$, where $L$ is the total length of the curve.  We first compute $L$ by evaluating 
\begin{equation}
	L = \int_0^{2\pi} \sqrt{\left(\frac{dx}{d\theta}\right)^2 + \left(\frac{dy}{d\theta}\right)^2} \, d\theta,
\end{equation}
where $dx/d\theta$ and $dy/d\theta$ are computed by differentiating (\ref{coords}), using the Fourier series representation of $f$ to compute its derivative.  This integral can be evaluated accurately using any high order quadrature desired.  We use Gauss-Legendre with 1000 nodes.  

One then sets $\theta_1 = 0$, and solves
\begin{equation}
	\frac{8N}{L} = \int_{\theta_{j-1}}^{\theta_j} \sqrt{\left(\frac{dx}{d\theta}\right)^2 + \left(\frac{dy}{d\theta}\right)^2} \, d\theta
\end{equation}
for $\theta_j$ given $\theta_{j-1}$ using Newton's method.  The integral is again evaluated using any desired quadrature.  We use Gauss-Legendre with 16 nodes.  

\end{appendix}

\bibliographystyle{plain}
\bibliography{AccDerivsPaper}

\end{document}